\documentclass[reqno,11pt]{amsart}
\usepackage{amsmath,amssymb,tabularx,setspace,color}
\usepackage{amsmath,amssymb,tabularx,setspace,color}
\usepackage{graphicx,graphics,enumerate}
\usepackage[margin=3cm]{geometry}
\usepackage{epstopdf}

\usepackage{graphicx}
\usepackage{caption}

\usepackage{subcaption}

\makeatletter
\renewcommand\section{\@startsection {section}{1}{\z@}%
                                   {-3.5ex \@plus -1ex \@minus -.2ex}%
                                   {2.3ex \@plus.2ex}%
                                   {\normalfont\large\bfseries}}
\makeatother
\theoremstyle{plain}
\newtheorem{theorem}{Theorem}[section]
\newtheorem{Example}{Example}[section]
\newtheorem{Theorem}{Theorem}

\newtheorem{Lemma}{Lemma}

\begin{document}
\doublespace
\title[ ]{Non-parametric estimation of conditional quantiles  for time series with heavy tails.}
\author[]%
{D\lowercase{eemat} C M\lowercase{athew}$^{\lowercase{a}}$,  H\lowercase{areesh}  G$^{\lowercase{b}}$ \lowercase{and}  S\lowercase{udheesh} K K\lowercase{attumannil}$^{\lowercase{c},\dag}$ \\   $^{\lowercase{a}}$S\lowercase{t.} T\lowercase{homas} C\lowercase{ollege}, P\lowercase{ala},  I\lowercase{ndia},\\
$^{\lowercase{b}}$N\lowercase{aval} P\lowercase{hysical} \lowercase{and} O\lowercase{ceanographic} L\lowercase{aboratory}, C\lowercase{ochin}, I\lowercase{ndia},\\$^{\lowercase{a}}$I\lowercase{ndian} S\lowercase{tatistical} I\lowercase{nstitute},
  C\lowercase{hennai}, I\lowercase{ndia}.}
%\thanks{$^{\dag}${Corresponding author E-mail: \tt skkattu@isichennai.res.in}.}

\begin{abstract}We propose a modified weighted Nadaraya-Watson estimator for the conditional distribution of a time series with heavy tails. We establish the asymptotic normality of the proposed estimator.  Simulation study is carried out to assess the performance of the estimator. We illustrate our method using a dataset.
\end{abstract}
\maketitle
\vspace{-0.3in}
\section{Introduction}\vspace{-0.1in}
Forecasting future observations is one of the most important problems in time series analysis. Estimating the conditional distribution or quantiles can solve prediction problems to some extent. Furthermore, the estimation of conditional quantiles has received particular attention in recent decades, finding applications in various fields such as econometrics, finance and related areas. One of the most common approaches in forecasting is to postulate certain parametric models. However, this may not yield good results due to the lack of complete information about the functional form of the model. Nonparametric approaches can overcome this problem by allowing different functional forms.

Nonparametric kernel-based smoothing provides estimators with desirable asymptotic properties. For comprehensive coverage of these techniques see Fan (2018) and Silverman (2018). In this context, the failure of least square-based methods has led researchers to search for more robust alternatives. Some of these robust methods have their roots in approaches suggested by Hardle (1984) and Hardle and Gasser (1984). Fan and Hall (1994) studied local median smoothing for independent data. Nonparametric estimation of conditional medians is discussed in Hall et al. (2002) and Chaudhuri and Loh (2002). For independent and identically distributed random variables, Stone (1977) established the weak consistency of the kernel estimates of the conditional median and Gannoun (1989)  proved its asymptotic normality.  Samanta (1989) extended the results for the estimates of the conditional quantiles. Under mixing assumptions, Boente and Fraiman (1995) proved the convergence of nonparametric estimates of the median.

For finding prediction intervals in time series, Hall et al. (1999) proposed a weighted Nadaraya-Watson estimator for estimating conditional distribution. This estimator modifies the Nadaraya-Watson estimator by introducing probabilities as weights, which satisfy certain constraints. The empirical likelihood method is applied to select these probabilities from the data. Under some regularity conditions, Cai (2002) established the asymptotic normality and weak consistency of these estimators for $\alpha$-mixing time series. The principle of maximum entropy provides an unbiased method for selecting the probability distribution when only partial information is available (Jaynes (1957), Cover and Thomas (1991)). This criterion is explored in this study to determine the weights used in estimating the conditional distribution and conditional quantile.

The rest of the paper is organized as follows: In Section 2, we discuss the estimation of conditional distribution and conditional quantile. We propose a modified weighted Nadaraya-Watson estimator for estimating the conditional distribution of a time series with heavy tails. The principle of maximum entropy is applied to select the optimal weights for estimating the conditional distribution in the proposed kernel regression. The estimators of conditional quantiles are used to predict future values of the time series and construct prediction intervals.  Section 3 establishes the weak consistency and asymptotic normality of these estimators under some regularity conditions. We also determine the optimal bandwidth of the kernel by minimizing the mean squared error.  Section 4 presents a Monte Carlo simulation study to assess the performance of the estimator. Finally, we illustrate the method using real data.
\vspace{-0.2in}
\section{Weighted Nadaraya-Watson estimator} \vspace{-0.1in}
Consider a pair of random variables $(Y,Z)$ defined on  $\mathbb{R}^p \times \mathbb{R}$. Consider a $\mathbb{R}^p \times \mathbb{R}$ valued strict stationary
process $\{(Y_i,Z_i)\}$ with distribution same as that of $(Y,Z)$.  The conditional distribution of $Z|Y$ is given by $F_{Z|Y}(z)=E(I(Z_i<z|Y_i=y))$ and we interpret this as regression of  $Z_i$ on $Y_i.$

For a given kernel function $K(.)$, let $p_i,\, 1\leq i \leq n $
(functions of $Y_i's$ ) be weights satisfying \vspace{-0.1in}
\begin{eqnarray}\label{eq6.1}
p_i \geq 0, \quad \sum_{i=1}^{n} p_i=1 \text{ and } \sum_{i=1}^n p_i (Y_i-y)
K_h(Y_i-y)=0.
\end{eqnarray}
The weighted  Nadayara-Watson estimator of the conditional distribution is given by (Hall et al., 1999) \vspace{-0.2in}
\begin{eqnarray}\label{eq6.2}
\widehat F(z|y)=\frac{\sum_{i=1}^n p_i
I(Z_i<z)K_h(Y_i-y)}{\sum_{i=1}^n p_i K_h(Y_i-y)}.
\end{eqnarray}
The third condition in (\ref{eq6.1}) requires weights that enforce local linearity in the kernel function. It is evident that $ \widehat
F(z|y)$ is a proper distribution function. Note that the weights
$p_i'$ s satisfying the above conditions are not uniquely determined. Hall et al. (1999) used the empirical likelihood method to obtain the weights. Here, we estimate the weights $p_i's$ based on the principle of maximum entropy.
%Then we study the asymptotic properties of the estimators of the conditional distributions and conditional quantile.

To find the weights $p_i'$s we maximize entropy subject to the constraints in (\ref{eq6.1}). For this purpose, consider the Lagrangian function given by \vspace{-0.1in}
\begin{eqnarray}\label{eq6.3}
G=-\sum_{i=1}^n p_i \log p_i + k \Big(\sum_{i=1}^{n} p_i-1\Big
)+\lambda\sum_{i=1}^n p_i (Y_i-y) K_h(Y_i-y).
\end{eqnarray}
Differentiating with respect to $k$ and $\lambda$ and equating to zero leads to the equations %\vspace{-0.1in}
\begin{eqnarray}\label{eq6.4}
&\sum_{i=1}^n  (Y_i-y) K_h(Y_i-y) e^{ \lambda (Y_i-y) K_h(Y_i-y) } =0 &\\\label{eq6.45}
&\sum_{i=1}^n e^{ \lambda (Y_i-y) K_h(Y_i-y) } = e^{1-k}.\vspace{-0.2in}
\end{eqnarray}
We obtain the weights as $p_i=e^{-1+k+\lambda (Y_i-y) K_h(Y_i-y)},$ where  $k$ and $\lambda$ can be computed from (\ref{eq6.4}) and (\ref{eq6.45}) using numerical methods.
\subsection{Conditional quantile estimators}
Consider a stationary process $\{(Y_i,Z_i)\}$ with the same distribution as that of $(Y,Z)$ on a
probability space $(\Omega, \mathbb{F},\mathbb{P}),$  where $Y_i$ is the lagged values of $Z_i.$ Our goal is to predict the future
values of $Z_{t+m}$ or to obtain the prediction interval for
$Z_{t+m}$ , for $m=1,2,\ldots$ from the past values
$Z_{t},Z_{t-1},\ldots$ .  The $\tau$-{th}  quantile can be obtained as the quantity which minimises, the pinball loss function (Steinwart and Christmann, 2011).
For the pinball loss function given by \vspace{-0.1in}
\begin{eqnarray}\label{eq6.6}
 L(Z,g(Y))=\tau(Z-g(Y))I(Z>g(Y))+(1-\tau)(g(Y)-Z)I(g(Y)\geq Z ),
\end{eqnarray}\vspace{-0.1in}
the $\tau$-{th}  quantile $q_{\tau}(Z|Y)$ satisfies
\begin{equation}\label{pinq}\vspace{-0.1in}
  q_{\tau}(Z|Y)=\mathrm{argmin}\{E(L(Z,g(Y)))\}.%\vspace{-0.1in}
\end{equation}This also motivates us to consider the prediction interval based on conditional quantiles.

For $\alpha \in (0,1)$ , a $100(1-\alpha)\%$ prediction interval can be constructed by taking  $\tau=\alpha/2$ and $\tau=1-\alpha/2$  in (\ref{pinq}). We then  obtain the prediction  interval as \vspace{-0.1in}
 $$[ q_{(\alpha/2)}(Z|Y), q_{(1-\alpha/2)}(Z|Y).$$
 Now, we have the estimator of $\tau$-{th} conditional
 quantile as an inverse of the proposed estimator of the conditional distribution.
 Let $F_{Z|Y}(.)$ be the conditional distribution function of $Z$ given
$Y.$ Then the $\tau$-{th}  conditional quantile is defined as \vspace{-0.1in}
\begin{equation}\label{eq6.5}
  q_{\tau}(Z|Y)=\inf\{x\in R : F_{Z|Y}(x)\geq \tau\}.\vspace{-0.1in}
\end{equation}
If $F_{Z|Y}$ is  strictly increasing we can write
$q_{\tau}(Z|Y)=F_{Z|Y}^{-1}(\tau)$. The $\tau$-{th} conditional quantile can be estimated by inverting the
 conditional distribution estimator given in (\ref{eq6.2}). Then for $\alpha \in (0,1)$, a $100(1-\alpha)\%$  prediction interval can be constructed from the estimated distribution function as $$ [\widehat
F_{Z|Y}^{-1}(\alpha/2),\widehat F_{Z|Y}^{-1}(1-\alpha/2)].$$
\vspace{-0.2in}
\section{Asymptotic Properties} \vspace{-0.1in}
Here, we discuss the asymptotic properties of the proposed estimators of conditional distribution and conditional quantile under some mixing conditions. Let $f(.)$ denote the marginal density of $Y_t.$ Define $k_j=\int u^j K(u) du,$  $v_0=\int K^2(u)du$, $ v_j= \int u^j K^2(u) du $ and
$\omega(z|y)^2=F(z|y)(1-F(z|y))/f(y).$ Let $s_j(y)=\frac{1}{nh^j}
\sum_{t=1}^n (Y_t-y)^j K_h(Y_t-y)$ and $S_n()$ denote the $m\times m$ matrix
with $s_{i+j-2}(y)$ as $(i,j)^{th}$ element. We assume the following regularity conditions.

\noindent A1. For fixed $z$ and $y$, $f(y)>0$ and $0<F(z|y)<1$, $f$ is
 continuous at $y$ and $F(y|.)$ has continuous second derivative in the neighborhood of $y$.

\noindent A2. The kernel $K$ is symmetric, compactly supported probability
density satisfying
\begin{equation*}
  |K(y_1)-K(y_2)|<C|y_1-y_2| \quad \text{for any } y_1 \text{ and } y_2.
\end{equation*}

\noindent A3. The process $(Y_t,Z_t)$ is regular in the sense
$$\beta(j)=\sup_{i>1}E\big\{\sup_{A\in \mathfrak{F}_{i+j}^\infty }P(A|\mathfrak{F}_1^i)-P(A)\big \} \rightarrow 0 \ \mathrm{as} \ j\rightarrow\infty, $$ where $\mathfrak{F}_i^j$ is the $\sigma$-field
generated by $\{(Y_k,Z_k) i \leq k \leq j\}.$
Also $\sum_{j\geq1} j^2
\beta(j)^{\delta/(\delta+1)}<\infty$ for some $\delta \in [0,1).$

\noindent A4. As $n \rightarrow \infty, h\rightarrow 0$ we have $\liminf_{n
\rightarrow \infty} nh^4>0$.

\noindent A5. For $t>1,$ let $f_{1,t}(.,.)$ be the joint density of $(Y_1,Y_t)$ and assume $|f_{1,t}(u,v)-f(u)f(v)|\leq C $ for all $u$ and $v.$\\
The proofs of the following theorem are given in  the  Appendix.
\begin{theorem}\label{thm6.1}\vspace{-0.2in}
Under the regularity conditions $ A1-A5 $, as $n \rightarrow \infty$
\begin{small}
\begin{eqnarray}\label{eq6.7}
\hat{F}(z|y)-F(z|y)=(nh)^{-1/2} \omega(z|y) v_0 ^{1/2} N + \frac
{1}{2} h^2 k_2 F''(z|y)+o_p\big(h^2+(nh)^{-1/2}\big),
\end{eqnarray}\vspace{-0.1in}
\end{small}where $N$ is a standard normal random variable.
Equivalently, we can state
\begin{eqnarray}\label{eq6.8}\vspace{-0.1in}
\sqrt{nh}\big [\hat{F}(z|y)-F(z|y)-\frac{1}{2} h^2 k_2 F''(z|y) \big
]
 \stackrel{\mathbb{L}}{\rightarrow} N(0,v_0 \omega ^2(z|y)).
\end{eqnarray}
\end{theorem}

Since the proposed conditional distribution estimator is monotone
increasing and lies between 0 and 1, using (\ref{eq6.5}) we can
obtain the estimator of  the $\tau$-{th} conditional quantile as
$\widehat{q}_\tau(Z|Y)$ satisfying
$\widehat{F}_{Z|Y}(\widehat{q}_\tau(Z|Y))=\tau,$ so that \vspace{-0.1in}
$$ \widehat{q}_\tau(Z|Y)=\inf\{z \in R : \widehat{F}_{Z|Y}(z) \geq \tau \}.\vspace{-0.1in}$$
%We can not use this inverting approach to
%`double kernel" local linear estimator due to lack of monotonicity
%and some other difficulties due to use of two bandwidths.

\noindent To prove the asymptotic properties of the estimator of the
conditional quantile, we assume the following additional conditions.

A6: Assume $f(z|y)$, the conditional density of $Z|Y$  exists and is
continuous at $y$.

A7: $f(q_\tau|y)>0$.\vspace{-0.2in}
\begin{theorem}\label{thm6.2}
Under the  regularity conditions $A1-A7$, as $n \rightarrow \infty$
\begin{eqnarray}\label{eq6.9}\vspace{-0.2in}
\hat q_\tau(z|y)\stackrel{P}{\rightarrow} q_\tau(z|y)\vspace{-0.2in}
\end{eqnarray}
and\vspace{-0.2in}
\begin{eqnarray*}\vspace{-0.2in}
\sqrt{nh}\big [\hat q_\tau(z|y)-q_\tau(z|y)-\frac{1}{2} \frac{h^2
k_2 F''(q_\tau |y)} {f(q_\tau |y)} \big ]\vspace{-0.2in}
 \stackrel{\mathbb{L}}{\rightarrow} N(0,\sigma^2(q_\tau)), \vspace{-0.2in}
\end{eqnarray*} where $\sigma^2(q_\tau) = \frac{v_0 \tau(1-\tau)}{f^2(q_\tau|y)g(q_\tau)}$.
\end{theorem}

  %Optimal bandwidth can be obtained using Theorem \ref{thm6.2}.
 Using  Theorem \ref{thm6.2}, we have the bias and the asymptotic variance of conditional quantile estimator $\widehat q_\tau$ as
$${Bias}^2(\widehat q_\tau)= \frac{1}{2} \frac{h^2 k_2 F''(q_\tau |y)}
{f(q_\tau |y)} $$ and $${Var}(\widehat q_\tau)=  \frac{\sigma^2(q_\tau)}{nh}
= \frac{1}{{nh}}\frac{v_0 \tau(1-\tau)}{f^2(q_\tau|y)g(q_\tau))} .$$ Therefore,  the
mean square error of $\widehat q_\tau$ is given by
\begin{equation}\label{eq6.10}
  {\rm{MSE}}( \widehat q_\tau ) = \Big [ \frac{h^2 k_2 F''(q_\tau |y)}
{2 f(q_\tau |y)} \Big]^2 + \frac{1}{{nh}} \frac{v_0
\tau(1-\tau)}{f^2(q_\tau|y)g(q_\tau)}.
\end{equation}
The optimal bandwidth  associated with conditional quantile can be obtained by minimizing (\ref{eq6.10}).
Differentiating (\ref{eq6.10}) with respect to $h$ and equating to
zero, we obtain the optimal bandwidth as
 \begin{equation*}
   h_{opt}=\Big[\frac{(k_2 F''(q_\tau |y))^2}{v_0 \tau(1-\tau)/g(q_\tau)}\Big]^{1/5}n^{-1/5}.
 \end{equation*}
\section{Simulation and data analysis}
The performance of the conditional distribution estimators is assessed through a numerical study and by analyzing real data.

For numerical illustration, we consider the AR(1) model $y_t=0.76y_{t-1}+\varepsilon_t$ with
$\varepsilon_t\sim N(0,1).$  We generate 500 observations from this model. The first 495 observations are used for estimating the conditional distribution and $95\%$ predictive interval for the last five observations obtained. The estimated conditional distribution of $Y_t|Y_{t-1}$ for a particular value is given in Figure 1(A).  In Table 1, we provide the predictive interval for the last five observations. It can be seen that all prediction intervals contain the true value.
\begin{figure}[h]
        \begin{subfigure}[b]{0.5\textwidth}
                \centering
                \includegraphics[width=6cm]{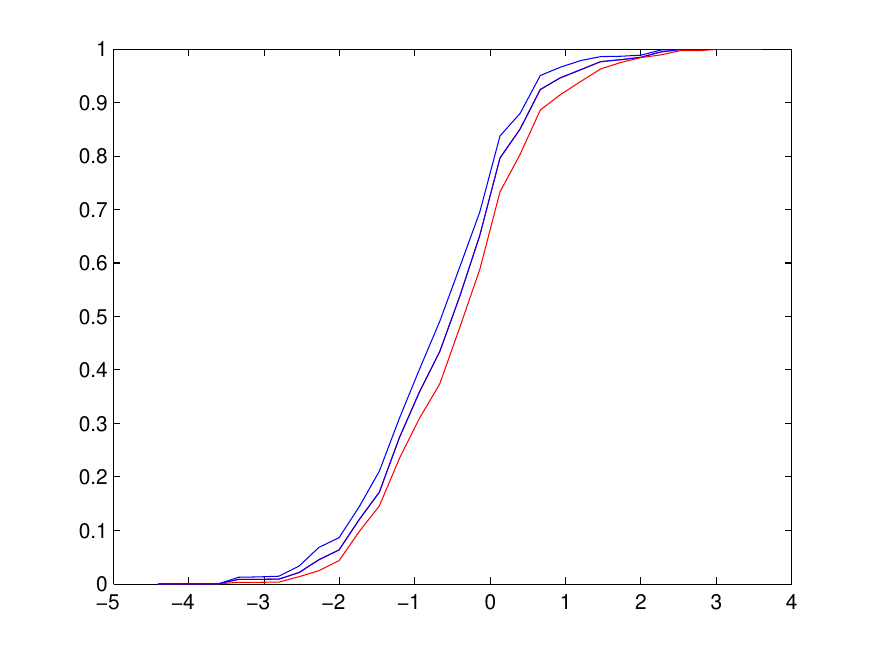}
                \caption{The conditional CDFs}
                \label{fig:1.1}
        \end{subfigure}%
         \caption{Conditional CDF estimated and predictive interval }\label{fig:1}
\end{figure}

\begin{table}[h]\footnotesize
\caption{Predictive interval for AR simulated data.}
\begin{tabular}{|r|r|r|r|}

\hline
           &     True Value &          Predictive Interval  \\
\hline
  $Y_{496}$ &      -1.002 &                    $[-2.20,2.08]$ \\
\hline
 $ Y_{497}$ &      0.2654&                     $[-2.40,1.92]$\\
\hline
 $Y_{498}$ &      0.0796&                     $[-2.36,1.02]$\\
\hline
  $Y_{499} $&      1.482 &                      $[-1.86,2.12]$ \\
\hline
   $Y_{500}$ &        0.8462 &                $[-2.26,1.90]$ \\
\hline
\end{tabular}
\end{table}
\begin{figure}[h]
        \begin{subfigure}[b]{0.5\textwidth}
                \centering
                \includegraphics[width=6cm]{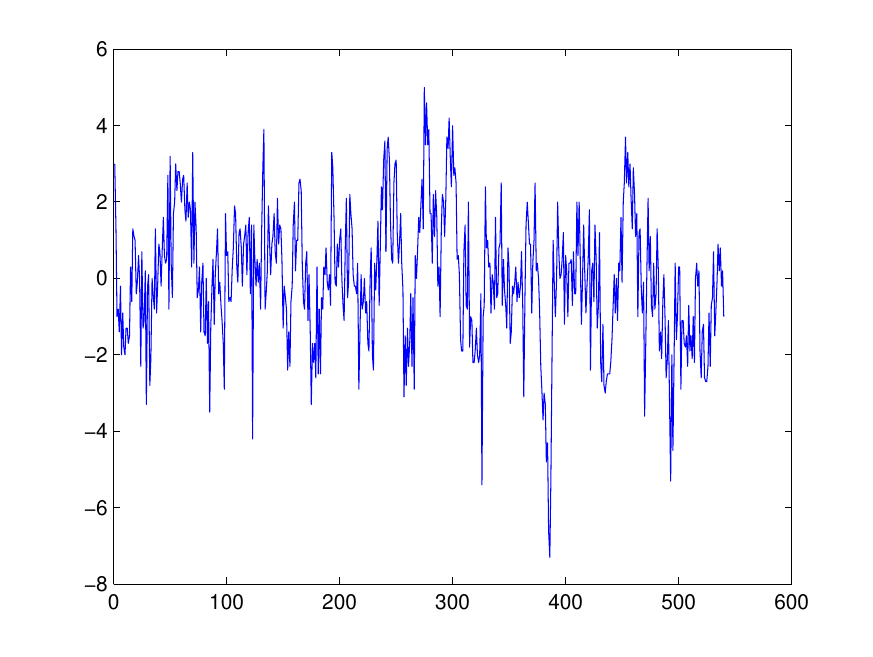}
                 \caption{Time Series Plot}
                \label{fig:1.1}
        \end{subfigure}%
        ~ %add desired spacing between images, e. g. ~, \quad, \qquad etc.
        \begin{subfigure}[b]{0.5\textwidth}
                \centering
                \includegraphics[width=6cm]{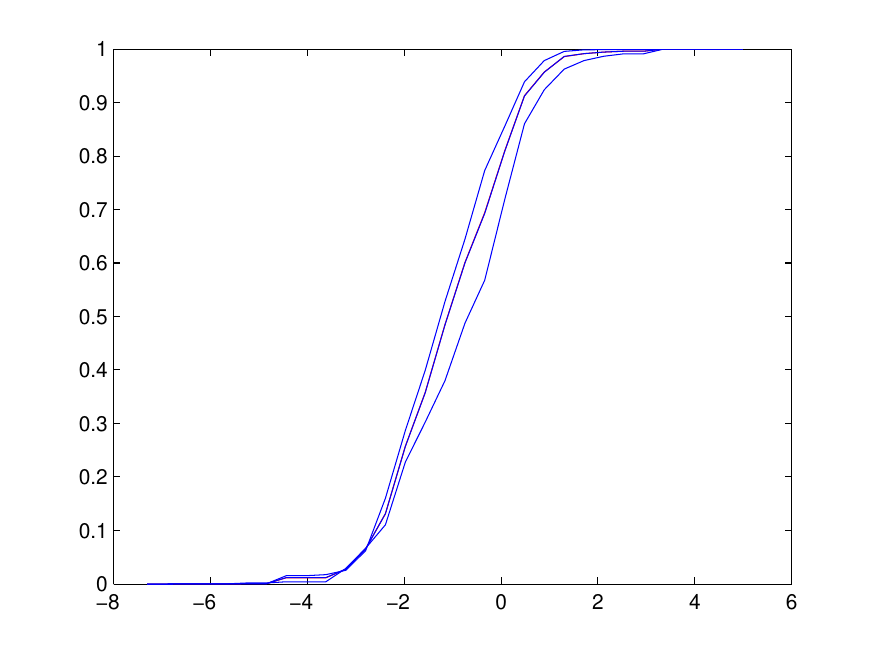}
                \caption{The conditional CDFs}
                \label{fig:1.2}
        \end{subfigure}
        \caption{Time Series Plot and estimated conditional CDF  }\label{fig:1}
\end{figure}

Next, we consider the monthly values of the Southern Oscillation Index (SOI) during 1950-1995. This series consists of 540 observations on the SOI, computed as the "difference of the departure from the long-term monthly mean sea level pressures" at Tahiti in the South Pacific and Darwin in Northern Australia. Figure 2(A) shows the time series plot of the data, and it can be assumed to be stationary. The first 535 observations are used for estimation, and the last 5 observations are left for prediction.  The conditional distribution
of $Y_t$ given $X_t=Y_{t-1}$ is shown in Figure 2(B) for three different values of $X_t.$  Now we consider the forecasting of the last five observations. Table 2 provides the $95\%$ prediction interval for the last 5 observations. It can be seen that in all cases, the prediction intervals contain the true values.

\begin{table}[h]\footnotesize
\caption{Predictive interval for SOI data.}
\begin{tabular}{|r|r|r|r|}

\hline
           &     True Value &           Predictive Interval  \\
\hline
  $Y_{536}$ &      0.19374 &                   $[-1.42,1.53]$ \\
\hline
 $ Y_{537}$ &      0.79374&                 $[-0.71,1.99]$\\
\hline
 $Y_{538}$ &      -0.20696&                   $[-1.69,2.13]$\\
\hline
  $Y_{539} $&      0.19374 &                     $[-1.56,1.62]$ \\
\hline
   $Y_{540}$ &        -1.006 &                 $[-2.61,0.59]$ \\
\hline
\end{tabular}
\end{table}

%\section{Conclusions}

\vspace{-0.1in}
\section*{Appendix}\vspace{-0.1in}
We need following Lemma to prove Theorem \ref{thm6.1} and the proof follows from the ergodic theorem (see Hall et al. 1999).
\begin{Lemma}\label{lem1}For a given kernel $K$, define $k_j=\int u^j K(u) du$.
Under the regularity conditions A1-A5, $S_n(y)\rightarrow f(y)S(y)$ in probability, where $S()$ denote the $m\times m$ matrix with $k_{i+j-2}$ as
$(i,j)^{th}$ element.
\end{Lemma}
%\newpage
\noindent \textbf{Proof of Theorem \ref{thm6.1}}: Consider
\begin{eqnarray}\label{decomp}
\widehat{F}_{Z|Y}-F_{Z|Y} &=& \frac{\sum_{i=1}^n p_i
I(Z_i<z)K_h(Y_i-y)}{\sum_{i=1}^n p_i K_h(Y_i-y)} - F_{Z|Y}\nonumber \\
&=& \frac{\sum_{i=1}^n[
I(Z_i<z)-F_{Z|Y}]p_i K_h(Y_i-y)}{\sum_{i=1}^n p_i K_h(Y_i-y)} \nonumber\\
&=&\big\{(nh)^{-1}T_1+T_2 \big\} T_3^{-1} \{1+o_p(1)\},
\end{eqnarray}
where
$$T_1=\sqrt{\frac{h}{n}}\sum_{t=1}^n e ^{-\lambda (Y_i-y) K_h(Y_i-y) } [I(Z_i<z)-F_{Z|Y_t}] K_h(Y_i-y),$$
$$T_2=\sum_{t=1}^n (F_{Z|Y_t}-F_{Z|Y})$$ and
$$T_3=\frac{1}{n}\sum_{t=1}^n e ^{-\lambda (Y_i-y) K_h(Y_i-y)}
K_h(Y_i-y). $$

To prove the theorem, first, using Lindeberg-Feller central limit theorem, we  show $T_1$ converges in distribution to normal random variable. We then prove that $T_2$
converges in probability to the bias term in Theorem \ref{thm6.1} and
$T_3$ converges in probability to a bounded quantity.

First, we find the mean and the variance of $T_1.$ Consider
\begin{eqnarray*}
T_1&=&\sqrt{\frac{h}{n}}\sum_{t=1}^n e ^{-\lambda (Y_i-y) K_h(Y_i-y) } \varepsilon_t K_h(Y_i-y)\\
&=&\frac{1}{\sqrt{n}}\sum_{t=1}^n \xi_t,
\end{eqnarray*}
where $\varepsilon_t=[I(Z_i<z)-F_{Z|Y_t}]$ and $\xi_t= \sqrt{h} e
^{-\lambda (Y_i-y) K_h(Y_i-y) } \varepsilon_t K_h(Y_i-y) $. It can be
easily verified that $E(\xi_t)=0$.
\begin{eqnarray*}
E(T_1^2)&=&\frac{1}{n}E\left(\sum_{i=1}^n \xi_i^2\right)+\frac{1}{n}\sum_{i=2}^n \left(1-\frac{i-1}{n}\right) Cov \left(\xi_1,\xi_i\right)\\
&=&E( \xi_1^2)+\frac{1}{n}\sum_{i=2}^n \left(1-\frac{i-1}{n}\right) Cov(\xi_1,\xi_i).
\end{eqnarray*}
\begin{eqnarray*}
E(\xi_t^2)&=&E(h e ^{-2\lambda (Y_i-y) K_h(Y_i-y) } \varepsilon_t^2 K_h^2(Y_i-y))\\
&=&E(E(h e ^{-2\lambda (Y_i-y) K_h(Y_i-y) } \varepsilon_t^2 K_h^2(Y_i-y)|Y_t))\\
&=&F(z|y)(1-F(z|y))E(h e ^{-2\lambda (Y_i-y) K_h(Y_i-y) }
K_h^2(Y_i-y)).
\end{eqnarray*}
Applying Lemma \ref{lem1} in the Taylor series expansion of $e^x$, we obtain
\begin{eqnarray}\label{eq9.26}
\nonumber E(\xi_t^2)&=& F(z|y)(1-F(z|y)) v_0 f(y) + o(1) \\
&=& v_0 \omega^2(z|y) f^2(y).
\end{eqnarray}
Choose $d_n=O(h^{-1/(1+\delta/2)})$ and write
\begin{eqnarray*}
 \sum_{i=2}^n \left(1-\frac{i-1}{n}\right) Cov (\xi_1,\xi_i)&=&\sum_{i=2}^{d_n} \left(1-\frac{i-1}{n}\right) Cov (\xi_1,\xi_i)\\&&+\sum_{i=d_n+1}^n \left(1-\frac{i-1}{n}\right) Cov (\xi_1,\xi_i).
\end{eqnarray*}
Applying theorem A.5 of Hall and Heyde (1980) and the assumption A2 on the
kernel, we obtain
$$|Cov(\xi_1,\xi_i)|\leq C h^{-1} \alpha(i-1),$$ which implies
$$ \sum_{i=d_n+1}^n (1-\frac{i-1}{n}) Cov (\xi_1,\xi_i) \leq C h^{-1} d_n^{-(1+\delta)}=o(1).$$
Due to the assumption A5, we have
$$ \sum_{i=2}^{d_n} \left(1-\frac{i-1}{n}\right) Cov (\xi_1,\xi_i)=O(d_n h)=o(1).$$
Hence, $E(T_1^2) \rightarrow v_0 \omega^2(z|y) f^2(y) $.

Next we prove the asymptotic normality of $T_1.$ Partition
$\{1,2,\ldots,n\}$ into $2q_n+1$ subsets with blocks of sizes
$r=\lfloor (nh)^{1/2} \rfloor$ and $s=\lfloor (nh)^{1/2}/\log n
\rfloor$ where $q=\lfloor \frac{n}{r_n+s_n}\rfloor$, $\lfloor x\rfloor$ denotes the greatest integer less than $x$.

Write $T_1$ as
\begin{eqnarray*}
T_1&=& \frac{1}{\sqrt n} \big[ \sum_{j=0}^{q-1}\eta_j + \sum_{j=0}^{q-1}\zeta_j +\eta_q \big] \\
&=&\frac{1}{\sqrt n} \big[ T_{1,1}+T_{1,2}+T_{1,3} \big],
\end{eqnarray*}
where $$\eta_j =\sum_{i=j(r+s)}^{j(r+s)+r-1} \xi_i,\,\,
\zeta_j=\sum_{i=j(r+s)+r}^{(j+1)(r+s)} \xi_i\,\text{ and }
\eta_q=\sum_{i=q(r+s)}^{n-1} \xi_i.$$

Consider  $ \frac{1}{\sqrt n}  T_{1,1} = \frac{1}{\sqrt n}
\sum_{j=0}^{q-1}\eta_j  $.  By Lemma 1.1 of Volkonskii and Rozanov
(1959), we obtain
\begin{eqnarray*}
 \Bigg| E \exp (itT_{1,1})-\prod_{j=0}^{q-1} E \exp (it\eta_j) \Bigg|  &\leq& 16 \frac n r \alpha(s) \\
  &\rightarrow& 0,\,\, \text {as } n\rightarrow \infty,
 \end{eqnarray*}%\vspace{-0.2in}
which shows that $\eta_j's$ in $T_{1,1}$ are asymptotically independent.

Similar arguments used in obtaining the variance of $T_1$,
we have
\begin{eqnarray*}
E(\frac{1}{ n} T_{1,1}^2) &=& \frac{1}{n}\sum_{j=0}^{q-1} E (\eta_j)^2 \\
&=& \frac{q}{n}Var(\sum_{j=1}^r \xi_j) \rightarrow v_0 \omega^2(z|y) f^2(y).
\end{eqnarray*}

Hence by mixing conditions and Theorem 4.1 of Shao and Yu (1996), we
have
\begin{eqnarray*}
E\left(\eta_1^2 I(|\eta_1|\geq\varepsilon \omega(z|y)\sqrt n)\right) &\leq& C n^{-1/2} E(|\eta_1|^3) \\
&\leq& C n^{-1/2}r^{3/2} \{E|\xi_1|^6\}^{1/2}  \\
&\leq&  C n^{-1/2}r^{3/2}h^{-1}.
\end{eqnarray*}
Therefore, as $n\rightarrow \infty$
$$\frac{1}{n}\sum_{j=0}^{q-1} E\big[ \eta_1^2 I(|\eta_1|\geq\varepsilon \omega(z|y)\sqrt n)\big]\leq C(nh^3)^{-1/4}\rightarrow 0$$

By Lindeberg-Feller central limit theorem, $\frac{1}{\sqrt n}
T_{1,1}$ has asymptotic normal distribution. Next, consider $T_{1,2}=\sum_{j=0}^{q-1} \zeta_j$ and we have
\begin{eqnarray*}
E(T_{1,2}^2) &=&\sum_{j=0}^{q-1}Var(\zeta_j)+2\sum_{j>i} Cov(\zeta_i,\zeta_j)\\
&=& T_{1,2,1}+T_{1,2,2} \quad(say).
\end{eqnarray*}\vspace{-0.1in}
Once can easily verify\vspace{-0.1in}
$$T_{1,2,1}=q Var(\zeta_1)=q Var\left(\sum_{j=1}^s \xi_j\right)\rightarrow q s v_0 \omega^2(z|y) f^2(y)\vspace{-0.1in}  $$
and $$\vspace{-0.1in}
T_{1,2,2}\leq2\sum_{i=1}^{n-r}\sum_{j=i+r}^{n}|Cov(\xi_i,\xi_j)|=o(n).
$$Hence $$ \frac{1}{n} E(T_{1,2}^2)= o(1). $$ By similar arguments as above, we have
 $$ \frac{1}{n} E(T_{1,3}^2)= o(1). $$
From these we can conclude that $T_1$ is asymptotically normally distributed.

Now, consider the second term in the decomposition (\ref{decomp})  \vspace{-0.1in}$$T_2=\sum_{t=1}^n
(F_{Z|Y_t}-F_{Z|Y} ) .\vspace{-0.1in}$$ \vspace{-0.1in} By Taylor's theorem, we obtain
\begin{eqnarray*}\vspace{-0.1in}
T_2&=&\sum_{i=1}^n (F_{Z|Y_i}-F_{Z|Y} )\\
&=& \frac{1}{2n} \sum_{i=1}^n F^{(2)}_{Z|Y}(y_t-y)^2e ^{-\lambda (Y_i-y) K_h(Y_i-y)} K_h(Y_i-y)  + o_p(h^2).
\end{eqnarray*}
Hence by Lemma 1 given in the Appendix, we have $$T_2 \xrightarrow{p} h^2 \mu_2\, F^{(2)}_{Z|Y}\, g(y).$$\vspace{-0.1in}
Finally, we obtain
\begin{eqnarray*}\vspace{-0.1in}
T_3&=&\frac{1}{n}\sum_{i=1}^n e ^{-\lambda (Y_i-y) K_h(Y_i-y)} K_h(Y_i-y) \\
&=&\frac{1}{n}\sum_{i=1}^n \left(1-\lambda (Y_i-y) K_h(Y_i-y)+\frac{1}{2}\lambda^2 (Y_i-y)^2 K^2_h(Y_i-y)+\ldots\right)K_h(Y_i-y)\\
&\rightarrow&g(y)+o_p(1).
\end{eqnarray*}
Combining all the above results, we have the proof of the theorem.

\noindent \textbf{Proof of Theorem \ref{thm6.2}}:
Since $\widehat{F}_Z|Y$ is monotone increasing, from Theorem
\ref{thm6.1} above and Theorem 1 of Tucker (1967),  as $n\rightarrow\infty$\vspace{-0.1in}
\begin{eqnarray}\vspace{-0.3in}
\sup_{y\in R} |\widehat{F}_Z|Y(z|y)-F_Z|Y(z|y)| \xrightarrow{p} 0.
\end{eqnarray}\vspace{-0.1in}
Due to the  uniqueness of the quantile, for the choice of
$ \delta=\min\{\tau-F(q_\tau-\varepsilon|y),F(q_\tau+\varepsilon|y)-\tau\}>0 ,$
we have \vspace{-0.2in}
\begin{eqnarray*}
P(|\widehat{q}_\tau-q_\tau|>\varepsilon)&\leq& P(|F(\widehat q_\tau|y)-\tau|>\delta)\\
&=& P(|\widehat F(\widehat q_\tau|y)- F(\widehat q_\tau|y|>\delta)\\
&\leq& P(\sup_{y\in R} |\widehat{F}_Z|Y(z|y)-F_Z|Y(z|y)|>\delta)
\rightarrow 0,
\end{eqnarray*}
and the asymptotic normality is immediate from Theorem \ref{thm6.1}.
\end{document}